\documentclass[12pt,twoside]{article}
\usepackage[english]{babel}
\usepackage[latin1]{inputenc}
\usepackage{amsmath}
\usepackage{amssymb,amsfonts}
\usepackage{graphicx}
\newcommand{\bL}{\mathbf{L}}

\newcommand{\bR}{\mathbf{R}}

\newcommand{\bS}{\mathbf{S}}

\newcommand{\bx}{\mathbf{x}}

\newcommand{\cS}{\mathcal{S}}

\newcommand{\HYP}{\mathbb{H}^3}
\newcommand{\HYN}{\mathbb{H}^n}

\newcommand{\SLR}{\widetilde{\bS\bL_2\bR}}
\newcommand{\NIL}{\mathbf{Nil}}
\newcommand{\SOL}{\mathbf{Sol}}

\begin{document}
\pagestyle{myheadings}
\markboth{\centerline{Jen\H o Szirmai}}
{Upper bound of density $\dots$}
\title
{Upper bound of density for packing of congruent hyperballs in hyperbolic $3-$space}

\author{\normalsize{Jen\H o Szirmai} \\
\normalsize Budapest University of Technology and \\
\normalsize Economics, Institute of Mathematics, \\
\normalsize Department of Geometry \\
\date{\normalsize{\today}}}

\maketitle

%%%%%%%%%%%%%%%%%%%%%%%%%%%%%%%%%%%%%%%%%%%%

\begin{abstract}
In \cite{Sz17-2} we proved that 
to each saturated congruent hyperball packing exists a decomposition of $3$-dimensional hyperbolic space $\HYP$ into truncated tetrahedra. 
Therefore, in order to get a density upper bound for hyperball packings, it is sufficient to determine
the density upper bound of hyperball packings in truncated simplices.
In this paper we prove, using the above results and the results of papers \cite{M94} and \cite{Sz14}, that the density upper bound of the saturated 
congruent hyperball (hypersphere) packings related to the
corresponding truncated tetrahedron cells is realized in a regular truncated tetrahedra with density $\approx  0.86338$.
Furthermore, we prove that the density of locally optimal congruent hyperball arrangement 
in regular truncated tetrahedron is not a monotonically increasing function of the height (radius) of corresponding optimal hyperball, 
contrary to the ball (sphere) and horoball (horosphere) packings.
\end{abstract}

%%%%%%%%%%%%%%%%%%%%%%%%%%%%%%%%%%%%%%%%%%%
\newtheorem{theorem}{Theorem}[section]
\newtheorem{corollary}[theorem]{Corollary}
\newtheorem{conjecture}{Conjecture}[section]
\newtheorem{lemma}[theorem]{Lemma}
\newtheorem{exmple}[theorem]{Example}
\newtheorem{defn}[theorem]{Definition}
\newtheorem{rmrk}[theorem]{Remark}
            %%% for no-italic, numbered environments, use:
\newenvironment{definition}{\begin{defn}\normalfont}{\end{defn}}
\newenvironment{remark}{\begin{rmrk}\normalfont}{\end{rmrk}}
\newenvironment{example}{\begin{exmple}\normalfont}{\end{exmple}}
            %%% for unnumbered environments, use f.e.
\newenvironment{acknowledgement}{Acknowledgement}
%%%%%%%%%%%%%%%%%%%%%%%%%%%%%%%%%%%%%%%%%%%%%%%%%%%%%%%%%%%%%%%%%%%%

%============================================================================%
%                             the main article                               %
%============================================================================%

%%%%%%%%%%%%%%%%%%%%%%%%%%%%%%%%%%%%%%%%%%%%%%%%%%%%%%%%%%%%%%%%%%%%%%%%%%%%%%
\section{Preliminary results}
Let $X$ denote a space of constant curvature, either the $n$-dimensional sphere $\mathbb{S}^n$,
Euclidean space $\mathbb{E}^n$, or hyperbolic space $\mathbb{H}^n$ with $n \geq 2$. An important question of discrete geometry
is to find the highest possible packing density in $X$ by congruent non-overlapping balls of a given radius \cite{Be}, \cite{G--K}.

Euclidean cases are the best explored. One major recent development has been the settling of the long-standing Kepler conjecture,
part of Hilbert's 18th problem, by Thomas Hales at the turn of the 21st century.
Hales' computer-assisted proof was largely based on a program set forth by L. Fejes T\'oth in the 1950's \cite{Ha}.

In $n$-dimensional hyperbolic geometry occur several new questions concerning the  packing and covering problems, e.g.
in $\HYN$ there are $3$ kinds of ``generalized balls (spheres)":  the {\it usual balls} (spheres),
{\it horoballs} (horospheres) and {\it hyperballs} (hyperspheres). Moreover, the definition of packing density is crucial in hyperbolic spaces
as shown by B\"or\"oczky \cite{B78}, for standard examples also see \cite{G--K}, \cite{R06}.
The most widely accepted notion of packing density considers the local densities of balls with respect to their Dirichlet--Voronoi cells
(cf. \cite{B78} and \cite{K98}). In order to consider ball packings in $\overline{\mathbb{H}}^n$, we use an extended notion of such local density.

In space $X^n$ let $d_n(r)$ be the density of $n+1$ mutually touching spheres or horospheres of radius
$r$ (in case of horosphere $r=\infty$) with respect
to the simplex spanned by their centres. L.~Fejes T\'oth and H.~S.~M.~Coxeter
conjectured that the packing density of balls of radius $r$ in $X^n$ cannot exceed $d_n(r)$.
This conjecture has been proved by C.~A.~Rogers for Euclidean space $\mathbb{E}^n$.
The 2-dimensional spherical case was settled by L.Fejes T\'oth in \cite{FTL}. 

\vspace{3mm}

{\bf Ball (sphere) and horoball (horosphere) packings:}

\vspace{3mm}

In \cite{B78} and \cite{B--F64} K.~B\"or\"oczky proved the following theorem
for {\it ball and horoball} packings for any $n$ ($2 \le n \in \mathbb{N})$:
\begin{theorem}[K.~B\"or\"oczky]
In an $n$-dimensional space of constant curvature consider a packing of spheres of radius $r$.
In spherical space suppose that $r<\frac{\pi}{4}$.
Then the density of each sphere in its Dirichlet-Voronoi cell cannot
exceed the density of $n+1$ spheres of radius $r$ mutually
touching one another with respect to the simplex spanned by their centers.
\end{theorem}

The above greatest density in $\mathbb{H}^3$ is $\approx 0.85328$
which is not realized by packing with equal balls. However, it is attained by the horoball packing
(in this case $r=\infty$) of
$\overline{\mathbb{H}}^3$ where the ideal centers of horoballs lie on the
absolute figure of $\overline{\mathbb{H}}^3$. This ideal regular
tetrahedron tiling is given with Coxeter-Schl\"afli symbol $\{3,3,6\}$.
Ball packings of hyperbolic $n$-space and of other Thurston geometries
are extensively discussed in the literature see e.g. \cite{Be}, \cite{B78},
\cite{G--K--K}, \cite{MoSzi18} and \cite{Sz14-1}, where the reader finds further references as well.

In a previous paper \cite{KSz} we proved that the above known optimal horoball packing
arrangement in $\mathbb{H}^3$ is not unique using the notions of horoballs in same and different types.
Two horoballs in a horoball packing are of the "same type" iff the
local densities of the horoballs to the corresponding cell
(e.g. D-V cell or ideal simplex) are equal, (see \cite{Sz12-2}).
We gave several new examples of horoball packing arrangements based on
totally asymptotic Coxeter tilings that yield the above B\"or\"oczky--Florian packing density upper bound (see \cite{B--F64})

We have also found that the B\"or\"oczky-Florian type density
upper bound for horoball packings of different types is no longer valid for fully
asymptotic simplices in higher dimensions $n > 3$  (see \cite{Sz12}).
For example in $\mathbb{H}^4$, the density of such optimal,
locally densest horoball packing is $\approx 0.77038$ larger than the
analogous B\"or\"oczky-Florian type density upper bound of $\approx 0.73046$.
However, these horoball packing configurations are only locally optimal and cannot be extended to the whole hyperbolic space $\mathbb{H}^4$.

In papers \cite{KSz14} and \cite{KSz18} we have continued our previous investigation in $\mathbb{H}^n$ ($n \in \{4,5\}$) allowing horoballs of different types.
We gave several new examples of horoball packing configurations that yield high densities ($\approx 0.71645$ in $\mathbb{H}^4$ and $\approx 0.59421$ in $\mathbb{H}^5$)
where horoballs are centered at ideal vertices of certain Coxeter simplices, and are invariant
under the actions of their respective Coxeter groups.

\vspace{3mm}

{\bf Hyperball (hypersphere) packings:}

\vspace{3mm}

A hypersphere is the set of all points in
$\HYN$, lying at a certain distance, called its {\it height}, from a hyperplane, on both sides of the hyperplane (cf. \cite{V79} for the planar case).

In hyperbolic plane $\mathbb{H}^2$ the universal upper bound of the hypercycle packing density is $\frac{3}{\pi}$,
and the universal lower bound of hypercycle covering density is $\frac{\sqrt{12}}{\pi}$,
proved by I.~Vermes in \cite{V72, V79, V81}. We note here that independently from him in \cite{MM03} T.~H. Marshall and G.~J. Martin obtained
similar results to the hypercycle packings.

In \cite{Sz06-1} and \cite{Sz06-2} we analysed the regular prism tilings (simply truncated Coxeter orthoscheme tilings) and the corresponding optimal hyperball packings in
$\mathbb{H}^n$ $(n=3,4)$ and we extended the method developed of paper \cite{Sz06-2} to
$5$-dimensional hyperbolic space (see \cite{Sz13-3}).
In paper \cite{Sz13-4} we studied the $n$-dimensional hyperbolic regular prism honeycombs
and the corresponding coverings by congruent hyperballs and we determined their least dense covering densities.
Furthermore, we formulated conjectures for candidates of the least dense hyperball
covering by congruent hyperballs in $3$- and $5$-dimensional hyperbolic spaces.

In \cite{Sz17-1} we discussed congruent and non-congruent hyperball packings of the truncated regular tetrahedron tilings.
These are derived from the Coxeter simplex tilings $\{p,3,3\}$ $(7\le p \in \mathbb{N})$ and $\{5,3,3,3,3\}$
in $3$- and $5$-dimensional hyperbolic space.
We determined the densest hyperball packing arrangement and its density
with congruent hyperballs in $\mathbb{H}^5$ and determined the smallest density upper bounds of
non-congruent hyperball packings generated by the above tilings in $\HYN,~ (n=3,5)$.

In \cite{Sz17} we deal with the packings derived by horo- and hyperballs (briefly hyp-hor packings) in $n$-dimensional hyperbolic spaces $\HYN$
($n=2,3$) which form a new class of the classical packing problems.
We constructed in the $2-$ and $3-$dimensional hyperbolic spaces hyp-hor packings that
are generated by complete Coxeter tilings of degree $1$
and we determined their densest packing configurations and their densities.
We proved using also numerical approximation methods that in the hyperbolic plane ($n=2$) the density of the above hyp-hor packings arbitrarily approximate
the universal upper bound of the hypercycle or horocycle packing density $\frac{3}{\pi}$ and
in $\HYP$ the optimal configuration belongs to the $\{7,3,6\}$ Coxeter tiling with density $\approx 0.83267$.
Furthermore, we analyzed the hyp-hor packings in
truncated orthosche\-mes $\{p,3,6\}$ $(6< p < 7, ~ p\in \mathbb{R})$ whose
density function is attained its maximum for a parameter which lies in the interval $[6.05,6.06]$
and the densities for parameters lying in this interval are larger that $\approx 0.85397$.

In \cite{Sz14} we proved that if the truncated tetrahedron is regular, then the density
of the densest packing is $\approx 0.86338$. This is larger than the B\"or\"oczky-Florian density upper bound
but our locally optimal hyperball packing configuration cannot be extended to the entirety of
$\mathbb{H}^3$. However, we described a hyperball packing construction,
by the regular truncated tetrahedron tiling under the extended Coxeter group $\{3, 3, 7\}$ with maximal density $\approx 0.82251$.

Recently, (to the best of author's knowledge) the candidates for the densest hyperball
(hypersphere) packings in the $3,4$ and $5$-dimensional hyperbolic space $\mathbb{H}^n$ are derived by the regular prism
tilings which have been studied in papers \cite{Sz06-1}, \cite{Sz06-2} and \cite{Sz13-3}.

In \cite{Sz17-2} we considered hyperball packings in
$3$-dimensional hyperbolic space and developed a decomposition algorithm that for each saturated hyperball packing provides a decomposition of $\HYP$
into truncated tetrahedra. Therefore, in order to get a density upper bound for hyperball packings, it is sufficient to determine
the density upper bound of hyperball packings in truncated simplices.

In \cite{Sz18} we studied hyperball packings related to truncated regular octahedron and cube tilings that are derived from the Coxeter simplex tilings
$\{p,3,4\}$ $(7\le p \in \mathbb{N})$ and $\{p,4,3\}$ $(5\le p \in \mathbb{N})$
in $3$-dimensional hyperbolic space $\HYP$. We determined the densest hyperball packing arrangement and its density
with congruent and non-congruent hyperballs related to the above tilings.
Moreover, we prove that the locally densest congruent or non-congruent hyperball configuration belongs to the regular truncated cube with density
$\approx 0.86145$. This is larger than the B\"or\"oczky-Florian density upper bound for balls and horoballs.
We described a non-congruent hyperball packing construction, by the regular cube tiling under the extended Coxeter group $\{4, 3, 7\}$
with maximal density $\approx 0.84931$.

In \cite{Sz18-1} we examined congruent and non-congruent hyperball packings generated by doubly truncated Coxeter orthoscheme tilings in the $3$-dimensional hyperbolic space.
We proved that the densest congruent hyperball packing belongs to the Coxeter orthoscheme tiling of parameter $\{7,3,7\}$ with density $\approx 0.81335$.
This density is equal -- by our conjecture -- with the upper bound density of the corresponding non-congruent hyperball arrangements.

\begin{rmrk}
If we try to define the density of system of sets in hyperbolic space as we did in Euclidean space,
i.e. by the limiting value of the density with respect to a sphere $C(r)$ of radius $r$ with a fixed centre $O$.
But since for a fixed value of $h$ the volume of spherical shell $C(r+h)-C(r)$ is the same order of magnitude as the volume of $C(r)$,
the argument used in Euclidean space to prove that the limiting value is independent of the choice of $O$ is does not work in hyperbolic space.
{\it Therefore the definition of packing density is crucial in hyperbolic spaces $\HYN$} as shown by K.~B\"or\"oczky \cite{B78}, for nice examples also see \cite{G--K}, \cite{R06}.
The most widely accepted notion of packing density considers the local densities of balls with respect to their Dirichlet--Voronoi cells
(cf. \cite{B78} and \cite{K98}), but in our cases these cells are infinite hyperbolic polyhedra.
The other possibility: the packing density $\delta$ can be defined (see \cite{V79}, \cite{V81}, \cite{Sz06-1},
\cite{Sz13-3}) as the reciprocal of the ratio of the volume of a fundamental
domain for the symmetry group of a tiling to the volume of the ball pieces contained in
the fundamental domain ($\delta<1)$. Similarly is defined the covering density $\Delta > 1$.
In the present paper our aim is to determine a density upper bound for saturated, congruent hyperball packings in $\HYP$ therefore we use an extended notion of such local density.
\end{rmrk}
\section{Saturated hyperball packings in $\HYP$ and their density upper bound}
We use for $\mathbb{H}^3$ (and analogously for $\HYN$, $n\ge3$) the projective model in the Lorentz space $\mathbb{E}^{1,3}$
that denotes the real vector space $\mathbf{V}^{4}$ equipped with the bilinear
form of signature $(1,3)$,
$
\langle \mathbf{x},~\mathbf{y} \rangle = -x^0y^0+x^1y^1+x^2y^2+ x^3 y^3,
$
where the non-zero vectors
$
\mathbf{x}=(x^0,x^1,x^2,x^3)\in\mathbf{V}^{4} \ \  \text{and} \ \ \mathbf{y}=(y^0,y^1,y^2,y^3)\in\mathbf{V}^{4},
$
are determined up to real factors, for representing points of $\mathcal{P}^n(\mathbb{R})$. Then $\mathbb{H}^3$ can be interpreted
as the interior of the quadric
$
Q=\{(\mathbf{x})\in\mathcal{P}^3 | \langle  \mathbf{x},~\mathbf{x} \rangle =0 \}=:\partial \mathbb{H}^3
$
in the real projective space $\mathcal{P}^3(\mathbf{V}^{4},
\mbox{\boldmath$V$}\!_{4})$ (here $\mbox{\boldmath$V$}\!_{4}$ is the dual space of $\mathbf{V}^{4}$).
Namely, for an interior point $\mathbf{y}$ holds $\langle  \mathbf{y},~\mathbf{y} \rangle <0$.

Points of the boundary $\partial \mathbb{H}^3 $ in $\mathcal{P}^3$
are called points at infinity, or at the absolute of $\mathbb{H}^3 $. Points lying outside $\partial \mathbb{H}^3 $
are said to be outer points of $\mathbb{H}^3 $ relative to $Q$. Let $(\mathbf{x}) \in \mathcal{P}^3$, a point
$(\mathbf{y}) \in \mathcal{P}^3$ is said to be conjugate to $(\mathbf{x})$ relative to $Q$ if
$\langle \mathbf{x},~\mathbf{y} \rangle =0$ holds. The set of all points which are conjugate to $(\mathbf{x})$
form a projective (polar) hyperplane
$
pol(\mathbf{x}):=\{(\mathbf{y})\in\mathcal{P}^3 | \langle \mathbf{x},~\mathbf{y} \rangle =0 \}.
$
Thus the quadric $Q$ induces a bijection
(linear polarity $\mathbf{V}^{4} \rightarrow
\mbox{\boldmath$V$}\!_{4})$
from the points of $\mathcal{P}^3$ onto their polar hyperplanes.

Point $X (\bold{x})$ and hyperplane $\alpha (\mbox{\boldmath$a$})$
are incident if $\bold{x}\mbox{\boldmath$a$}=0$ ($\bold{x} \in \bold{V}^{4} \setminus \{\mathbf{0}\}, \ \mbox{\boldmath$a$}
\in \mbox{\boldmath$V$}_{4}
\setminus \{\mbox{\boldmath$0$}\}$).

{\it The hypersphere (or equidistance surface)} is a quadratic surface at a constant distance
from a plane (base plane) in both halfspaces. The infinite body of the hypersphere, containing the base plane, is called hyperball.

\vspace{3mm}

The {\it half hyperball } with distance $h$ to a base plane $\beta$ is denoted by $\mathcal{H}^h_+$.
The volume of a bounded hyperball piece $\mathcal{H}^h_+(\mathcal{A})$,
delimited by a $2$-polygon $\mathcal{A} \subset \beta$, and its prism orthogonal to $\beta$, can be determined by the classical formula
(2.1) of J.~Bolyai \cite{B31}.
\begin{equation}
\mathrm{Vol}(\mathcal{H}^h_+(\mathcal{A}))=\frac{1}{4}\mathrm{Area}(\mathcal{A})\left[k \sinh \frac{2h}{k}+
2 h \right], \tag{2.1}
\end{equation}
The constant $k =\sqrt{\frac{-1}{K}}$ is the natural length unit in
$\mathbb{H}^3$, where $K$ denotes the constant negative sectional curvature. In the following we may assume that $k=1$.

Let $\mathcal{B}^h$ be a hyperball packing in $\HYP$ with congruent hyperballs of height $h$.

The notion of {\it saturated packing} follows from that fact that the density of any packing can be improved by adding further packing elements as long as there is
sufficient room to do so. However, we usually apply this notion for packings with congruent elements.

In \cite{Sz17-2} {\it we modified the classical definition of saturated packing} for non-compact ball packings with generalized balls 
(horoballs, hyperballs) in $n$-dimensional hyperbolic space $\HYN$ ($n\ge 2$ integer parameter):
\begin{defn}
A ball packing with non-compact generalized balls (horoballs or/and hyperballs) in $\HYN$ is saturated if no new non-compact generalized ball can be added to it.
\end{defn}
To obtain hyperball (hypersphere) packing upper bound it obviously suffices to study saturated hyperball packings (using the above definition) and in what follows we assume that
all packings are saturated unless otherwise stated.

We take the set of hyperballs $\{ \mathcal{H}^h_i\}$ of a saturated hyperball packing $\mathcal{B}^h$ (see Definition 2.1).
Their base planes are denoted by $\beta_i$.
Thus in a saturated hyperball packing the distance between two ultraparallel base planes
$d(\beta_i,\beta_j)$ is at least $2h$ (where for the natural indices holds $i < j$
and $d$ is the hyperbolic distance function).

In \cite{Sz17-2} we described a procedure to get a decomposition of 3-dimensional hyperbolic space $\HYP$ into truncated
tetrahedra corresponding to a given saturated hyperball packing whose main steps were the following:
\begin{enumerate}
\item Using the radical planes of the hyperballs $\mathcal{H}^h_i$, similarly to the Euclidean space, can be constructed
the unique Dirichlet-Voronoi (in short D-V) decomposition of $\HYP$ to the given hyperball packing $\mathcal{B}^h$. 

\item We consider an arbitrary {\it proper} vertex $P \in \HYP$ of the above $D-V$ decomposition and the hyperballs $\mathcal{H}^h_i(P)$
whose D-V cells
meet at $P$. The base planes of the hyperballs $\mathcal{H}^h_i(P)$ are denoted by $\beta_i(P)$, and these planes determine a
non-compact polyhedron $\mathcal{D}^i(P)$ with the intersection of their halfspaces
containing the vertex $P$. Moreover, denote $A_1,A_2,A_3,\dots$ the outer vertices of $\mathcal{D}^i(P)$ and cut off
$\mathcal{D}^i(P)$ with the polar planes $\alpha_j(P)$ of its outer vertices $A_j$. Thus, we obtain a convex compact polyhedron
$\mathcal{D}(P)$.
This is bounded by the base planes $\beta_i(P)$ and "polar planes" $\alpha_j(P)$. Applying this procedure for all vertices of the
above Dirichlet-Voronoi decomposition, we obtain an other decomposition of $\HYP$ into convex polyhedra.
\item We consider $\mathcal{D}(P)$ as a tile of the above decomposition. The planes from the finite set of base planes $\{\beta_i(P)\}$ are
called adjacent if there is a vertex $A_s$ of $\mathcal{D}^i(P)$ that lies on each of the above plane.
We consider non-adjacent planes $\beta_{k_1}(P),\beta_{k_{2}}(P),\beta_{k_{3}}(P), \dots \beta_{k_{m}}(P) \in \{\beta_i(P)\}$ $(k_l  \in \mathbb{N}^+, ~ l=1,2,3,\dots m)$
that have an outer point of intersection denoted by $A_{k_1\dots k_m}$. Let $N_{\mathcal{D}(P)} \in \mathbb{N}$ denote the {\it finite} number of the outer points $A_{k_1\dots k_m}$
related to $\mathcal{D}(P)$. It is clear, that its
minimum is $0$ if $\mathcal{D}^i(P)$ is tetrahedron.
The polar plane $\alpha_{k_1\dots k_m}$ of $A_{k_1\dots k_m}$ is orthogonal to planes $\beta_{k_1}(P),\beta_{k_2}(P), \dots \beta_{k_m}(P)$
(thus it contain their poles $B_{k_1}$, $B_{k_2}$,\dots $B_{k_m}$) and divides $\mathcal{D}(P)$ into two convex polyhedra
$\mathcal{D}_1(P)$ and $\mathcal{D}_2(P)$.
\item If $N_{\mathcal{D}_1(P)} \ne 0$ and $N_{\mathcal{D}_2(P)} \ne 0$ then $N_{\mathcal{D}_1(P)} < N_{\mathcal{D}(P)}$ and $N_{\mathcal{D}_2(P)} < N_{\mathcal{D}(P)}$ then
we apply the point 3 for polyhedra $\mathcal{D}_i(P),~ i \in \{1,2\}$.
\item If $N_{\mathcal{D}_i(P)} \ne 0$ or $N_{\mathcal{D}_j(P)} = 0$ ($i \ne j,~ i,j\in\{1,2\}$) then we consider the polyhedron $\mathcal{D}_i(P)$ where
$N_{\mathcal{D}_i(P)}=N_{\mathcal{D}(P)}-1$ because the vertex $A_{k_1\dots k_m}$ is left out and apply the point 3.
\item If $N_{\mathcal{D}_1(P)}=0$ and $N_{\mathcal{D}_2(P)}=0$ then the procedure is over for $\mathcal{D}(P)$. We continue the procedure with the next cell.
\item We have seen in steps 3, 4, 5 and 6 that the number of the outer vertices $A_{k_1 \dots k_m}$ of any polyhedron obtained 
after the cutting process is less than the original one, and
we have proven in step 7 that the original hyperballs form packings in the new polyhedra $\mathcal{D}_1(P)$ and $\mathcal{D}_2(P)$, as well.
We continue the cutting procedure described in step 3 for both polyhedra $\mathcal{D}_1(P)$ and $\mathcal{D}_2(P)$. 
If a derived polyhedron is a truncated
tetrahedron then the cutting procedure does not give new polyhedra, thus the procedure will not be continued.
Finally, after a {\it finite number of cuttings} we get a decomposition of $\mathcal{D}(P)$ into truncated tetrahedra,
and in any truncated tetrahedron the corresponding congruent hyperballs from $\{ \mathcal{H}^h_i\}$ form a packing. 
Moreover, we apply the above method for the further cells.
\end{enumerate}
From the above algorithm we obtained the following
\begin{theorem}[J.~Sz.~\cite{Sz17-2}]
The in \cite{Sz17-2} described algorithm provides for each congruent saturated hyperball packing a decomposition of $\HYP$ into truncated tetrahedra. ~ ~$\square$
\end{theorem}
\begin{rmrk}
Przeworski,~A. proved a similar theorem in \cite{P13} but it was true only for cases if the base planes of hyperspheres form ``symmetric cocompacts arrangements" in $\HYN$.
\end{rmrk}
In \cite{M94} Y.~Miyamoto proved the analogue theorem of K.~B\"or\"oczky's theorem (Theorem 1.1):
\begin{theorem}[Y.~Miyamoto, \cite{M94}]
If a region in $\HYN$ bounded by hyperplanes has a hyperball (hypersphere) packing of height (radius) $r$ about its boundary, then in some sense, 
the ratio of its volume to the volume of its boundary is at least that of a regular truncated simplex of (inner) edgelength $2r$.
\end{theorem}
\begin{rmrk}
Independently from the above paper A.~Przeworski proved a similar theorem with other methods in \cite{P13}.
\end{rmrk}
Therefore, in order to get density upper bound related to the saturated hyperball packings it is sufficient to determine
the density upper bound of hyperball packings in truncated regular simplices (see Fig.~1).
\begin{figure}[ht]
\centering
\includegraphics[width=8cm]{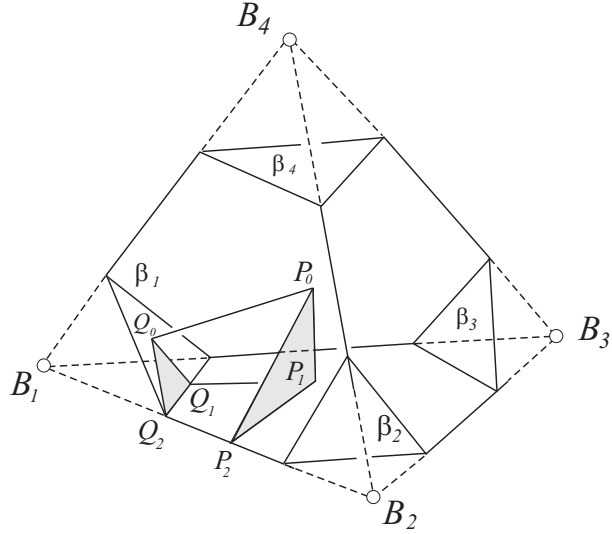} 
\caption{Regular truncated simplex, $\mathcal{S}(p)$, $p\in (6,\infty)$ with a simply truncated orthoscheme $\mathcal{O}=Q_0Q_1Q_2P_0P_1P_2$}
\label{}
\end{figure}
Thus, in the following we assume that the ultraparallel base planes $\beta_i$ of
$\mathcal{H}^{h(p)}_i$ $(i=1,2,3,4$, and $6<p\in \mathbb{R})$
generate a ``regular truncated tetrahedron" $\mathcal{S}(p)$ with outer vertices $B_i$
(see Fig.~1) whose non-orthogonal dihedral angles are equal to $\frac{2\pi}{p}$,
and the distances between two base planes $d(\beta_i,\beta_j)=:e_{ij}$  ($i < j \in \{1,2,3,4\})$ are equal to $2h(p)$ depending on the angle
$\frac{\pi}{p}$.

The truncated regular tetrahedron $\mathcal{S}(p)$ can be decomposed into $24$ congruent simply truncated orthoschemes; one of them
$\mathcal{O}=Q_0Q_1Q_2P_0P_1P_2$ is illustrated in Fig.~1 where $P_0$ is the centre of the ``regular tetrahedron" $\mathcal{S}(p)$,
$P_1$ is the centre of a hexagonal face of $\mathcal{S}(p)$, $P_2$ is the midpoint of a ``common perpendicular" edge of this face,
$Q_0$ is the centre of an adjacent regular triangle face of $\mathcal{S}(p)$, $Q_1$ is the midpoint of an appropriate edge of this face and
one of its endpoints is $Q_2$.

In our case the essential dihedral angles of orthoschemes $\mathcal{O}$ are the following:
$\alpha_{01}=\frac{\pi}{p}, \ \ \alpha_{12}=\frac{\pi}{3}, \ \
\alpha_{23}=\frac{\pi}{3}$. Therefore, the volume $\mathrm{Vol}(\mathcal{O})$ of the orthoscheme $\mathcal{O}$ and the volume
$\mathrm{Vol}(\mathcal{S}(p))=24 \cdot \mathrm{Vol}(\mathcal{O})$ can be computed for any given parameter $p$ $(6<p\in \mathbb{R})$ by Theorem 2.6.
\begin{theorem}[R.~Kellerhals, \cite{K89}, Theorem II.] The volume of a three-di\-men\-si\-o\-nal hyperbolic
complete orthoscheme (except Lambert cube cases, i.e. complete orthoschemes of degree $m=2$ with outer edge) $\mathcal{O} \subset \mathbb{H}^3$
is expressed with the essential
angles $\alpha_{01},\alpha_{12},\alpha_{23}, \ (0 \le \alpha_{ij} \le \frac{\pi}{2})$
in the following form:

\begin{align}
&\mathrm{Vol}(\mathcal{O})=\frac{1}{4} \{ \mathcal{L}(\alpha_{01}+\theta)-
\mathcal{L}(\alpha_{01}-\theta)+\mathcal{L}(\frac{\pi}{2}+\alpha_{12}-\theta)+ \notag \\
&+\mathcal{L}(\frac{\pi}{2}-\alpha_{12}-\theta)+\mathcal{L}(\alpha_{23}+\theta)-
\mathcal{L}(\alpha_{23}-\theta)+2\mathcal{L}(\frac{\pi}{2}-\theta) \}, \notag
\end{align}
where $\theta \in [0,\frac{\pi}{2})$ is defined by:
$$
\tan(\theta)=\frac{\sqrt{ \cos^2{\alpha_{12}}-\sin^2{\alpha_{01}} \sin^2{\alpha_{23}
}}} {\cos{\alpha_{01}}\cos{\alpha_{23}}},
$$
and where $\mathcal{L}(x):=-\int\limits_0^x \log \vert {2\sin{t}} \vert dt$ \ denotes the
Lobachevsky function.
\end{theorem}
In this case for a given parameter $p$ the length of the common perpendiculars $h(p)=\frac{1}{2}e_{ij}$ $(i < j$, $i,j \in \{1,2,3,4\})$
can be determined by the machinery of projective metric geometry. (In the following $\bx \sim c\cdot \bx$ with $c \in \mathbb{R}\setminus \{\mathbf{0}\}$
represent the same point $X = (\bx \sim c \cdot \bx)$ of $\mathcal{P}^3$.)

The points $P_2({\mathbf{p}}_2)$ and $Q_2({\mathbf{q}}_2)$ are proper points of hyperbolic $3$-space and
$Q_2$ lies on the polar hyperplane $pol(B_1)(\mbox{\boldmath$b$}^1)$ of the outer point $B_1$ thus

The hyperbolic distance $h(p)$ can be calculated by the following formula (see\cite{Sz14}):
\[
\begin{gathered}
\cosh{h(p)}=\cosh{P_2Q_2}=\frac{- \langle {\mathbf{q}}_2, {\mathbf{p}}_2 \rangle }
{\sqrt{\langle {\mathbf{q}}_2, {\mathbf{q}}_2 \rangle \langle {\mathbf{p}}_2, {\mathbf{p}}_2 \rangle}}= \\ =\frac{h_{23}^2-h_{22}h_{33}}
{\sqrt{h_{22}\langle \mathbf{q}_2, \mathbf{q}_2 \rangle}} =
\sqrt{\frac{h_{22}~h_{33}-h_{23}^2}
{h_{22}~h_{33}}},
\end{gathered} \notag
\]
where $h_{ij}$ is the inverse of the Coxeter-Schl\"afli matrix
\[(c^{ij}):=\begin{pmatrix}
1& -\cos{\frac{\pi}{p}}& 0&0\\
-\cos{\frac{\pi}{p}} & 1 & -\cos{\frac{\pi}{3}}&0\\
0 & -\cos{\frac{\pi}{3}} & 1&-\cos{\frac{\pi}{3}} \\
0&0&-\cos{\frac{\pi}{3}}&1\\
\end{pmatrix} \tag{2.2}
\]
of the orthoscheme $\mathcal{O}$.
We get that the volume $\mathrm{Vol}(\mathcal{S}(p))$, the maximal height $h(p)$ of the congruent hyperballs lying in $\mathcal{S}(p)$ and
$\mathrm{Vol}(\mathcal{H}^{h(p)} \cap \mathcal{S}(p))$ all depend only on
the parameter $p$ of the truncated regular tetrahedron $\mathcal{S}(p)$.

Therefore, the locally optimal density of congruent hyperball packing related to the regular truncated tetrahedron of parameter $p$ is
\begin{equation}
\delta(\mathcal{S}(p)):=\frac{4\cdot \mathrm{Vol}(\mathcal{H}^{h(p)} \cap \mathcal{S}(p))}{\mathrm{Vol}({\mathcal{S}(p)})}, \notag
\end{equation}
and $\delta(\mathcal{S}(p))$ depends only on
$p$ $(6<p\in \mathbb{R})$. Moreover,
the total volume of the parts of the four hyperballs lying in $\cS(p)$ can be computed
by formula (2.1), and the volume of
$\mathcal{S}(p)$ can be determined by Theorem 2.6.
\begin{figure}[ht]
\centering
\includegraphics[width=6cm]{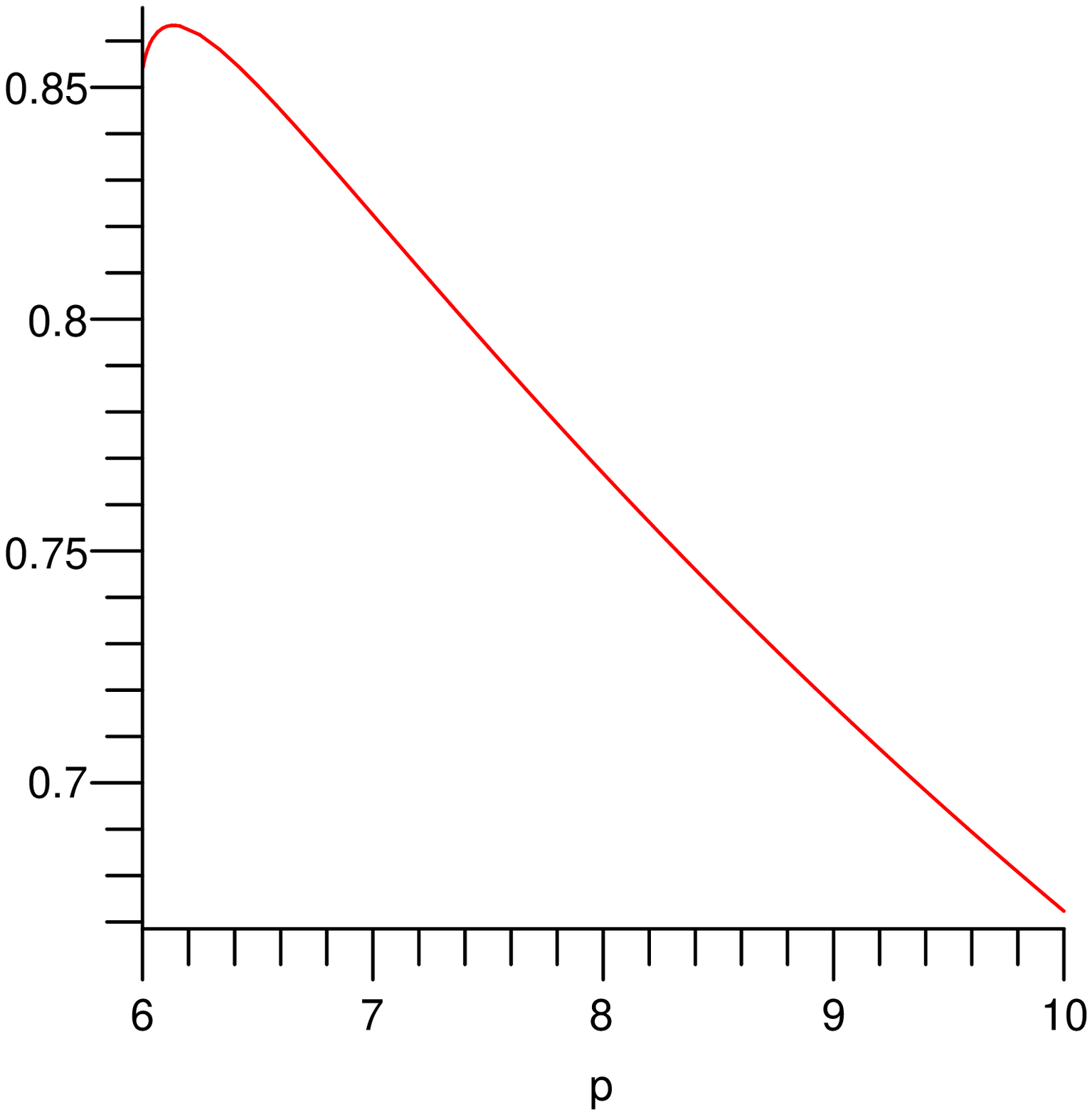} \includegraphics[width=6cm]{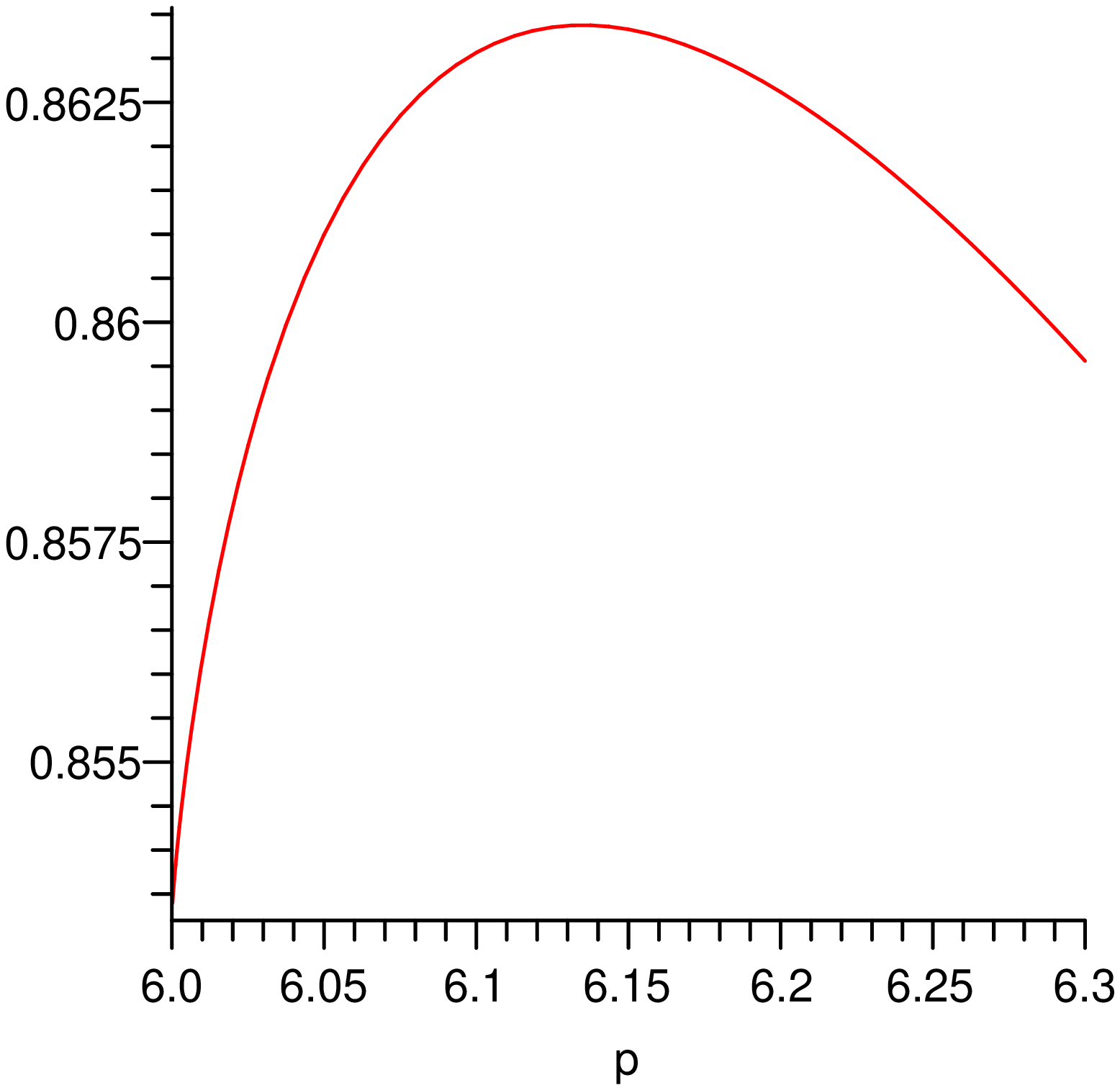}
\caption{The density function $\delta(\mathcal{S}(p))$, $p\in (6,10)$}
\label{}
\end{figure}
Finally, we obtain the plot after careful analysis of the smooth
density function (cf. Fig.~2) and we obtain the following
\begin{theorem}[J.~Sz.~\cite{Sz14}]
The density function $\delta(\mathcal{S}(p))$, $p\in (6,\infty)$
attains its maximum at $p^{opt} \approx 6.13499$, and $\delta(\mathcal{S}(p))$
is strictly increasing in the interval $(6,p^{opt})$, and strictly decreasing in $(p^{opt},\infty)$. Moreover, the optimal density
$\delta^{opt}(\mathcal{S}(p^{opt})) \approx 0.86338$ (see Fig.~2).
\end{theorem}
\begin{rmrk}
\begin{enumerate}
\item In our case $\lim_{p\rightarrow 6}(\delta(\mathcal{S}(p)))$ is equal to the B\"oröczky-Florian
upper bound of the ball and horoball packings in $\HYP$ \cite{B--F64} (observe that the dihedral angles of
$\cS(p)$ for the case of the horoball equal $2\pi/6$).
\item $\delta^{opt}(\mathcal{S}(p^{opt})) \approx 0.86338$ is larger than
the B\"oröczky-Florian upper bound
$\delta_{BF} \approx 0.85328$; but these hyperball packing configurations
are only locally optimal and cannot be extended to the entire hyperbolic
space $\mathbb{H}^3$.
\end{enumerate}
\end{rmrk}
We obtain the next theorem as the direct consequence of the previous statements:
\begin{theorem}
The density upper bound of the saturated congruent hyperball packings related to the
corresponding truncated tetrahedron cells is realized in a regular truncated tetrahedra belonging to parameter $p^{opt} \approx 6.13499$ with density $\approx  0.86338$.
\end{theorem}
We get from the above theorem directly the denial of the A.~Przeworski's conjecture \cite{P13}:
\begin{corollary}
The density function $\delta(\mathcal{S}(p))$, is not an increasing function of $h(p)$ (the height of hyperballs).
\end{corollary}
\begin{rmrk}
The hyperball packings in the regular truncated tetrahedra under the extended reflection groups with Coxeter-Schläfli symbol
$\{3,3,p\}$, investigated in paper \cite{Sz14}, can be extended to the entire hyperbolic space
if $p$ is an integer parameter bigger than $6$.
They coincide with the hyperball packings given by the regular $p$-gonal prism tilings in $\HYP$ with extended
Coxeter-Schl\"afli symbols $\{p,3,3\}$, see in \cite{Sz06-1}.
As we know, $\{3,3,p\}$ and $\{p,3,3\}$ are dually isomorphic extended reflection groups,
just with the above frustum of orthoscheme
as fundamental domain (Fig.~1, matrix $(c^{ij})$ in formula (2.2)).

In \cite{Sz14} we studied these tilings and the corresponding hyperball packings. Moreover, we computed their metric data for some 
integer parameters $p$ ($6<p\in \mathbb{N}$), where $\mathcal{A}$
is a trigonal face of the regular truncated tetrahedron, cf. Fig.~1. In the Table 1 we recalled from \cite{Sz14} important metric data of some ``realizable hyperball packings".
\medbreak
{\scriptsize
\centerline{\vbox{
\halign{\strut\vrule~\hfil $#$ \hfil~\vrule
&\quad \hfil $#$ \hfil~\vrule
&\quad \hfil $#$ \hfil\quad\vrule
&\quad \hfil $#$ \hfil\quad\vrule
&\quad \hfil $#$ \hfil\quad\vrule
\cr
\noalign{\hrule}
\multispan5{\strut\vrule\hfill\bf Table 1, \hfill\vrule}%
\cr
\noalign{\hrule}
\noalign{\vskip2pt}
\noalign{\hrule}
p & h(p) & \mathrm{Vol}(\mathcal{O}) & \mathrm{Vol}(\mathcal{H}^h_+(\mathcal{A}))& \delta(\mathcal{S}(p)) \cr
\noalign{\hrule}
7 & 0.78871 & 0.08856 & 0.07284 & 0.82251 \cr
\noalign{\hrule}
8 & 0.56419 & 0.10721 & 0.08220 & 0.76673 \cr
\noalign{\hrule}
9 & 0.45320 & 0.11825 & 0.08474 & 0.71663 \cr
\noalign{\hrule}
\vdots & \vdots  & \vdots  & \vdots  & \vdots \cr
\noalign{\hrule}
20 & 0.16397 & 0.14636 & 0.06064 & 0.41431 \cr
\noalign{\hrule}
\vdots & \vdots  & \vdots  & \vdots  & \vdots \cr
\noalign{\hrule}
50 & 0.06325 & 0.15167 & 0.02918 & 0.19240 \cr
\noalign{\hrule}
\vdots & \vdots  & \vdots  & \vdots  & \vdots \cr
\noalign{\hrule}
100 & 0.03147 & 0.15241 & 0.01549 & 0.10165 \cr
\noalign{\hrule}
p \to \infty & 0 & 0.15266 & 0 & 0 \cr
\noalign{\hrule}}}}}
\end{rmrk}
\medbreak
In hyperbolic spaces $\HYN$ ($n \ge 3$) the  problems of the densest ball, horoball and hyperball packings have not been settled yet,
in general (see e.g. \cite{KSz14}, \cite{Sz12}, \cite{Sz12-2}).
Moreover, the optimal sphere packing problem can be extended to the other homogeneous Thurston geometries, e.g. $\NIL$, $\SOL$, $\SLR$.
For these non-Euclidean geometries only very few results are known (e.g. \cite{Sz14-1} and the references given there).

By the above these we can say that the revisited Kepler problem keeps yet several interesting open questions.

I thank Prof. Emil Moln\'ar for his helpful comments and suggestions to this paper.
%============================================================================%
%                                references                                  %
%============================================================================%

\noindent
\footnotesize{Budapest University of Technology and Economics, Institute of Mathematics, \\
Department of Geometry, \\
H-1521 Budapest, Hungary. \\
E-mail:~szirmai@math.bme.hu \\
http://www.math.bme.hu/ $^\sim$szirmai}

\end{document}